\documentclass[11pt]{article}
\usepackage{amscd}
\usepackage{amsfonts}
\usepackage{amsmath}
\usepackage{amssymb}
\usepackage{amsthm}
\usepackage{bbm}
\usepackage{CJK}
\usepackage{fancyhdr}
\usepackage{graphicx}
\usepackage{hyperref}
\usepackage{indentfirst}
\usepackage{latexsym}
\usepackage{mathrsfs}
\usepackage{xypic}
\usepackage[compress]{cite}

\newtheorem{theorem}{Theorem}[section]

\newtheorem{definition}[theorem]{Definition}
\newtheorem{proposition}[theorem]{Proposition}
\newtheorem{example}[theorem]{Example}

\newtheorem{cor}[theorem]{Corollary}
\newtheorem{remark}[theorem]{Remark}

\usepackage[top=1in,bottom=1in,left=1.25in,right=1.25in]{geometry}
\def\<{\langle}
\def\>{\rangle}

\date{}
\begin{document}
\renewcommand{\baselinestretch}{1.2}
\renewcommand{\arraystretch}{1.0}
\title{\bf Splitting of operations for di-associative   algebras and  tri-associative algebras }
\author{{\bf Wen Teng}\\
{\small  School of Mathematics and Statistics, Guizhou University of Finance and Economics} \\
{\small  Guiyang  550025, P. R. of China}\\
  {\small E-mail: tengwen@mail.gufe.edu.cn } }
 \maketitle
\begin{center}
\begin{minipage}{13.cm}

{\bf Abstract:}
Loday  introduced di-associative   algebras and   tri-associative algebras  motivated  by  periodicity phenomena in algebraic $K$-theory.
   The purpose of this paper is to study the splittings of operations of di-associative   algebras and   tri-associative algebras.  First, we introduce the notion  of    a  quadri-dendriform algebra, which is a splitting of a di-associative algebra.
We show that   a relative averaging operator on dendriform algebras gives rise to a quadri-dendriform algebra.
 Conversely,  a quadri-dendriform algebra
gives rise to a  dendriform algebra and a  representation such that the quotient map is a relative averaging operator.
  Furthermore,    any quadri-dendriform algebra can be embedded into an averaging dendriform algebra.
Finally, we introduce the notion  of six-dendriform algebras,  which are a splitting of tri-associative algebras, and
 demonstrate that   homomorphic  relative averaging operators   induce   six-dendriform algebras.
 \smallskip

 \smallskip

 \end{minipage}
 \end{center}
 \normalsize\vskip0.5cm

\section{Introduction}
\def\theequation{\arabic{section}. \arabic{equation}}
\setcounter{equation} {0}

In \cite{Loday2001}, Loday introduced the concept of a di-associative algebra, inspired by his research on periodicity phenomena in algebraic $K$-theory.
The Koszul dual  of the operad associated with di-associative algebras is provided by the operad of dendriform algebras. A dendriform
algebra can be thought as a splitting of an associative algebra.
Recently, Das\cite{Das11802}  introduced the cohomology and deformation of dendriform algebras.
In \cite{Loday2004}, Loday and Ronco demonstrated that the family of chain modules over the standard simplices can be endowed with an operad structure. Algebras over this operad are referred to as tri-associative algebras, and they are distinguished by three   operations and eleven defining identities.
Additionally, they   introduced tridendriform algebras and established that tri-associative algebras are Koszul dual to tridendriform algebras.

The purpose of this paper is to study the splittings of operations for di-associative   algebras and tri-associative   algebras.
More precisely, we introduce the notion of quadri-dendriform algebras as a splitting  of   di-associative  algebras.
Subsequently, we introduce the notion of six-dendriform algebras as a splitting of   tri-associative  algebras.
The main tool to make connections between all these
structures is the so called averaging operator and more general   relative averaging operator.
The concept of the averaging operator was initially implicitly addressed by  Reynolds in 1895 \cite{Reynolds}.
In the last century, averaging operators   were mostly studied in the algebra of functions on a space  \cite{Kelley,Miller}.
In his PhD thesis,  Cao \cite{Cao}  examined averaging operators from both algebraic and combinatorial perspectives.
The further algebraic investigation of averaging operators on any binary operad, along with their relationships with bisuccessors, duplicators, and Rota-Baxter operators, is explicitly detailed in  \cite{Bai2013,Pei2020}.

The main results of this paper are summarized in the following diagram:
   $$\aligned
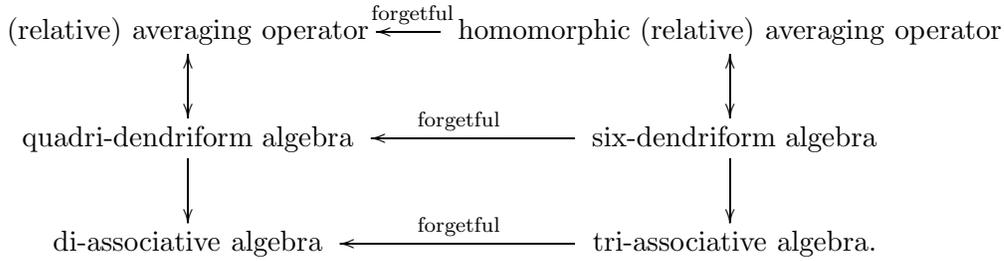
\xymatrix{
  \text{(relative) averaging operator} \ar[d] &\ar[l]_{\tiny{\text{forgetful}~~~~~~~~~~~~~~}} \text{ homomorphic  (relative) averaging operator } \ar[d]^{ } \\
  \ar[u] \text{ quadri-dendriform algebra }    \ar[d]  &\ar[l]_{\text{forgetful}} \ar[u]\text{  six-dendriform algebra}  \ar[d]\\
  \text{ di-associative   algebra }  & \ar[l]_{\text{forgetful}} \text{ tri-associative   algebra.}
   }
 \endaligned$$

  This paper is organized as follows.
In Section \ref{section:Quadri-dendriform algebras}, we introduce the notions
of  quadri-dendriform algebras and  relative averaging operators on  dendriform algebras. We demonstrate that  a relative averaging operator induces a quadri-dendriform algebra, and that  any
quadri-dendriform algebra can be embedded into an averaging  dendriform algebra.
In Section \ref{section:Six-dendriform algebras},  we  present the notion  of six-dendriform algebras, which are a splitting of   tri-associative algebras.
Additionally, we demonstrate that  homomorphic  relative averaging operators    give rise to   six-dendriform algebras.

Throughout this paper, $\mathbb{K}$ denotes a field of characteristic zero. All the  vector spaces  and
   (multi)linear maps are taken over $\mathbb{K}$.

\section{Quadri-dendriform algebras}\label{section:Quadri-dendriform algebras}
\def\theequation{\arabic{section}.\arabic{equation}}
\setcounter{equation} {0}
In this section,  we first introduce the notion of a quadri-dendriform algebra, which constitutes a splitting of a di-associative algebra. Next,
  we introduce   relative averaging operators on  dendriform algebras. We show that  a relative averaging operator induces a quadri-dendriform algebra. Finally, we show that any
quadri-dendriform algebra can be embedded into an averaging dendriform algebra.

\begin{definition}\cite{Loday2001}
A dendriform algebra is a vector space $D$  equipped  with two   binary operations  $\prec, \succ: D\otimes D\rightarrow D$  such that  for all $x, y, z\in D$,
\begin{align}
(x\prec y)\prec  z&=x\prec  (y\prec z+y\succ z), \label{2.1}\\
(x\succ  y)\prec  z&=x\succ  (y\prec z), \label{2.2}\\
x\succ (y\succ  z)&=(x\prec  y+x\succ  y)\succ z.\label{2.3}
\end{align}
\end{definition}
A morphism between two dendriform algebras is a linear map that preserves the corresponding
binary operations.   We denote by $\mathbf{Dend}$ the category
of dendriform algebras.

Dendriform algebras arise naturally from Rota-Baxter operator on associative algebras.
Let $(A,\mu)$ be an associative algebra and $R$ be a Rota-Baxter operator of it. In other words,
$R:A\rightarrow A$ is a linear map satisfying
$$\mu(Ra,Rb)=R(\mu(a,Rb)+\mu(Ra,b)),$$
for all $a,b\in A.$  In \cite{Aguiar}, Aguiar showed that if $R$ is a Rota-Baxter operator on $(A,\mu)$, then the two operations $\prec$ and $ \succ$ defined on $A$ by
$$a\prec b=\mu(a,Rb), ~~a\succ b=\mu(Ra,b),$$
defines a dendriform structure on $A$.

\begin{definition}\cite{Das11802}
Let $(D,\prec,\succ)$ be dendriform algebra.
A representation of $D$ is given by a vector space $V$ together with  four  bilinear maps
 \begin{align*}
 &\prec_{l}:A\otimes V\rightarrow V,~~~\succ_{l}: A\otimes V\rightarrow V,~~~\prec_{r}: V\otimes A\rightarrow V,~~~\succ_{r}: V\otimes A\rightarrow V,
 \end{align*}
 satisfying the following 9 identities
\begin{align}
&(x\prec y)\prec_{l} u=x\prec_{l}  (y\prec_{l} u+y\succ_{l} u),\label{2.4}\\
& (x\succ  y)\prec_{l}  u=x\succ_{l}  (y\prec_{l} u),\label{2.5}\\
& x\succ_{l} (y\succ_{l}  u)=(x\prec  y+x\succ  y)\succ_{l} u,\label{2.6}\\
&(x\prec_{l} u)\prec_{r}  z=x\prec_{l}  (u\prec_{r} z+u\succ_{r} z), \label{2.7}\\
&(x\succ_{l} u)\prec_{r}  z=x\succ_{l}  (u\prec_{r} z), \label{2.8}\\
&x\succ_{l} (u\succ_{r}  z)=(x\prec_{l}  u+x\succ_{l}  u)\succ_{r} z,\label{2.9}\\
&(u\prec_{r} y)\prec_{r}  z=u\prec_{r}  (y\prec z+y\succ z), \label{2.10}\\
&(u\succ_{r}  y)\prec_{r}  z=u\succ_{r}  (y\prec z), \label{2.11}\\
&u\succ_{r} (y\succ  z)=(u\prec_{r}  y+u\succ_{r}  y)\succ_{r} z,\label{2.12}
\end{align}
for all  $x,y,z\in D, u\in V$.
\end{definition}

Any dendriform algebra $(D,\prec,\succ)$  is a representation of itself with $\prec_{l}=\prec_{r}=\prec,\succ_{l}=\succ_{r}=\succ$.  This  representation is known as the adjoint  representation.

\begin{proposition}\label{proposition:  dendriform algebra}\cite{Das11802}
Let  $(D,\prec,\succ)$  be a dendriform algebra,  $V$ be a vector space, and let
$\prec_{l}:A\otimes V\rightarrow V, \succ_{l}: A\otimes V\rightarrow V $, $\prec_{r}: V\otimes A\rightarrow V $  and $ \succ_{r}: V\otimes A\rightarrow V$ be bilinear maps.
 Then $(V;\prec_{l},\succ_{l},\prec_{r},\succ_{r})$ is a representation
of $D$ if and only if $D\ltimes V : =(D\oplus V, \prec_\ltimes, \succ_\ltimes)$ is a  dendriform algebra, where  $ \prec_\ltimes $ and $ \succ_\ltimes$  are defined as:
 \begin{align*}
&(x,u)\prec_\ltimes(y, v):=(x\prec y, x\prec_l v+u\prec_r y),~&(x,u)\succ_\ltimes(y, v):=(x\succ y, x\succ_l u+v\succ_r y),
\end{align*}
for all  $x,y\in D, u, v\in V$. $D\ltimes V$ is called the semi-direct product of   dendriform algebra  $A$ with  $V$.
\end{proposition}

 In \cite{Loday2001} Loday introduced the notion of   di-associative algebras   in the study of
Leibniz algebras.

\begin{definition}\cite{Loday2001}
A di-associative algebra is a vector space $A$ equipped  with two   binary operations  $\vdash$ and $\dashv: A\otimes A\rightarrow A$,  such that  for all $x, y, z\in A$,
\begin{align*}
&(x\dashv y)\dashv  z=x\dashv  (y\vdash z), ~~~(x\dashv  y)\dashv  z=x\dashv  (y\dashv z),~~~(x\vdash  y)\dashv  z=x\vdash (y\dashv z),  \\
&(x\dashv  y)\vdash  z=x\vdash  (y\vdash z), ~~~(x\vdash  y)\vdash  z=x\vdash  (y\vdash z).
\end{align*}
\end{definition}

Di-associative algebras arise naturally from averaging operators on associative algebras.
Let $(A,\mu)$ be an associative algebra and $H$ be an averaging operator of it. In other words,
$H:A\rightarrow A$ is a linear map satisfying
$$\mu(Ha,Hb)=H\mu(a,Hb)=H\mu(Ha,b),$$
for all $a,b\in A.$  In \cite{Aguiar}, Aguiar showed that if $H$ is an averaging operator  on $(A,\mu)$, then the two operations $\vdash$ and $ \dashv$ defined on $A$ by
$$a\dashv b=\mu(a,Hb), ~~a\vdash b=\mu(Ha,b),$$
defines a di-associative algebra structure on $A$.

Inspired by Loday's\cite{Loday2001} work,
we now introduce a new algebraic structure:   quadri-dendriform algebras.

\begin{definition}
A quadri-dendriform algebra is a vector space $\mathcal{D}$ equipped with  four  binary operations $ \prec_\vdash,\prec_\dashv,\succ_\vdash,\succ_\dashv: \mathcal{D}\times\mathcal{D}\rightarrow\mathcal{D}$ such that
for all  $ x, y, z \in \mathcal{D}$,
\begin{align}
&(x\prec_\vdash y)\prec_{\vdash} z=(x\prec_\dashv y)\prec_{\vdash} z=x\prec_{\vdash}  (y\prec_{\vdash} z+y\succ_{\vdash} z),\label{2.13}\\
& (x\succ_\vdash y)\prec_{\vdash}  z=(x\succ_\dashv  y)\prec_{\vdash}  z=x\succ_{\vdash}  (y\prec_{\vdash} z),\label{2.14}\\
& x\succ_{\vdash} (y\succ_{\vdash}  z)=
\left\{ \begin{array}{ll}
=(x\prec_\vdash  y+x\succ_\vdash  y)\succ_{\vdash} z=(x\prec_\dashv  y+x\succ_\dashv  y)\succ_{\vdash} z,\\
=(x\prec_\vdash  y+x\succ_\dashv  y)\succ_{\vdash} z=(x\prec_\dashv  y+x\succ_\vdash  y)\succ_{\vdash} z,
 \end{array}  \right.\label{2.15}\\
&(x\prec_{\vdash} y)\prec_{\dashv}  z=x\prec_{\vdash}  (y\prec_{\dashv} z+y\succ_{\dashv} z), \label{2.16}\\
&(x\succ_{\vdash} y)\prec_{\dashv}  z=x\succ_{\vdash}  (y\prec_{\dashv} z), \label{2.17}\\
&x\succ_{\vdash} (y\succ_{\dashv}  z)=(x\prec_{\vdash}  y+x\succ_{\vdash}  y)\succ_{\dashv} z,\label{2.18}\\
&(x\prec_{\dashv} y)\prec_{\dashv}  z=
\left\{ \begin{array}{ll}
=x\prec_{\dashv}  (y\prec_\vdash z+y\succ_\vdash z)=x\prec_{\dashv}  (y\prec_\dashv z+y\succ_\dashv z),\\
=x\prec_{\dashv}  (y\prec_\vdash z+y\succ_\dashv z)=x\prec_{\dashv}  (y\prec_\dashv z+y\succ_\vdash z),
 \end{array}  \right.\label{2.19}\\
&(x\succ_{\dashv}  y)\prec_{\dashv}  z=x\succ_{\dashv}  (y\prec_\vdash z)=x\succ_{\dashv}  (y\prec_\dashv z), \label{2.20}\\
&x\succ_{\dashv} (y\succ_\vdash  z)=x\succ_{\dashv} (y\succ_\dashv z)=(x\prec_{\dashv}  y+x\succ_{\dashv}  y)\succ_{\dashv} z. \label{2.21}
\end{align}
\end{definition}

A morphism between two quadri-dendriform algebras is a linear map that
is compatible with the four  operations. We denote by  $\mathbf{QDend}$  the category
of quadri-dendriform algebras.


It is well known that for a given    dendriform algebra $(D,\prec, \succ)$, one can define a new binary operation $\diamond$ on $D$ by
$x\diamond y=x\prec y+x\succ y$ for $x,y\in D$. Consequently,   $(D,\diamond)$ forms an associative algebra.
Similarly,  the following conclusion can be drawn.

\begin{proposition}\label{proposition: di-associative algebra}
Let $(\mathcal{D},\prec_\vdash,\prec_\dashv,\succ_\vdash,\succ_\dashv)$   be a  quadri-dendriform algebra.
Then $(\mathcal{D}, \vdash, \dashv)$  is a  di-associative algebra, where
 \begin{eqnarray*}
&&x \vdash y=x\prec_\vdash y+x\succ_\vdash y,~~x\dashv y=x \prec_\dashv y+x\succ_\dashv y,
\end{eqnarray*}
for all $x,y\in  \mathcal{D}$.
\end{proposition}

\begin{proof}
For any $x,y,z\in  \mathcal{D}$, by Eqs.  \eqref{2.13}-\eqref{2.21}, we have
\begin{align*}
& (x\dashv y)\dashv  z=(x \prec_\dashv y+x\succ_\dashv y)\prec_\dashv z+(x \prec_\dashv y+x\succ_\dashv y)\succ_\dashv z\\
&=\left\{ \begin{array}{ll}
x\prec_{\dashv}  (y\prec_\dashv z+y\succ_\dashv z)+x\succ_{\dashv}  (y\prec_\dashv z)+x\succ_{\dashv} (y\succ_\dashv z)=x\dashv(y\dashv z),\\
x\prec_{\dashv}  (y\prec_\vdash z+y\succ_\vdash z)+x\succ_{\dashv}  (y\prec_\vdash z)+x\succ_{\dashv} (y\succ_\vdash  z)=x\dashv  (y\vdash z).
 \end{array}  \right.
\end{align*}
Similarly, we get
\begin{align*}
&(x\vdash  y)\dashv  z=(x\prec_\vdash y+x\succ_\vdash y)\prec_\dashv  z+(x\prec_\vdash y+x\succ_\vdash y)\succ_\dashv  z\\
&=x\prec_{\vdash}  (y\prec_{\dashv} z+y\succ_{\dashv} z)+x\succ_{\vdash}  (y\prec_{\dashv} z)+x\succ_{\vdash} (y\succ_{\dashv}  z)\\
&=x\vdash (y\dashv z)
\end{align*}
and
\begin{align*}
&(x\dashv  y)\vdash  z=(x \prec_\dashv y+x\succ_\dashv y)\prec_\vdash  z+(x \prec_\dashv y+x\succ_\dashv y)\succ_\vdash  z\\
&=x\prec_{\vdash}  (y\prec_{\vdash} z+y\succ_{\vdash} z)+x\succ_{\vdash}  (y\prec_{\vdash} z)+x\succ_{\vdash} (y\succ_{\vdash}  z)\\
&=x\vdash  (y\vdash z).
\end{align*}
Finally, we we obtain
\begin{align*}
&(x\vdash  y)\vdash  z=(x\prec_\vdash y+x\succ_\vdash y)\prec_\vdash  z+(x\prec_\vdash y+x\succ_\vdash y)\succ_\vdash  z\\
&=x\prec_{\vdash}  (y\prec_{\vdash} z+y\succ_{\vdash} z)+x\succ_{\vdash}  (y\prec_{\vdash} z)+x\succ_{\vdash} (y\succ_{\vdash}  z)\\
&=x\vdash  (y\vdash z).
\end{align*}
This completes the proof.
\end{proof}

Thus, according to Proposition \ref{proposition: di-associative algebra}, a  quadri-dendriform algebra can
be thought as a  di-associative algebra whose  two   operations are divided into four operations, and whose
di-associativities splits into nine novel identities. We denote by $\mathbf{Diass}$ the category of di-associative algebras. By the preceding proposition, there is a
well-defined functor:
$$\mathbf{QDend}\rightarrow\mathbf{Diass}.$$

\begin{example}
If $\prec_\vdash=\prec_\dashv=\prec$ and $\succ_\vdash=\succ_\dashv=\succ$, then we get simply a dendriform algebra.
So we get a functor between the categories of algebras:
$$\mathbf{Dend}\rightarrow\mathbf{QDend}.$$
\end{example}

\begin{example} \label{exam:3-tri-Leibniz kernel}
 Let $(\mathcal{D},\prec_\vdash,\prec_\dashv,\succ_\vdash,\succ_\dashv)$   be a  quadri-dendriform algebra.
The vector subspace $I_{\mathcal{D}}$ spanned by
\begin{align*}
\big\{x\prec_\vdash y-x\prec_\dashv y,~~x\succ_\vdash y-x\succ_\dashv y~:~x,y\in \mathcal{D}\big\}
\end{align*}
is an ideal of $\mathcal{D}$. It is clear that $\mathcal{D}$ is a dendriform algebra if and only if $I_{\mathcal{D}}=\{0\}$.
Therefore, the quotient algebra $\frac{\mathcal{D}}{I_{\mathcal{D}}}$ is a dendriform algebra.
\end{example}

\begin{example}  \label{exam:differential  dendriform algebra}
A differential dendriform algebra is a dendriform algebra $(D,\prec,\succ)$  equipped with a linear map
$d: D\rightarrow D$ satisfying
$$ d(x\prec y)= d(x)\prec y + x\prec d(y),~~d(x\succ y)= d(x)\succ y + x\succ d(y)  ~~\text{and}~~  d^2=0$$
for all $ x,y\in D$. Then it is easy to verify that $(\mathcal{D},\prec_\vdash,\prec_\dashv,\succ_\vdash,\succ_\dashv)$   is a  quadri-dendriform algebra, where
 $x\prec_\vdash y =d(x)\prec y, ~ x\prec_\dashv y= x\prec d(y),~ x\succ_\vdash y= d(x)\succ y$ and $x\succ_\dashv y=  x\succ d(y)$ for all $x,y\in D$.
\end{example}

\begin{example}
Let $(D,\prec,\succ)$  be a dendriform algebra and $(V;\prec_{l},\succ_{l},\prec_{r},\succ_{r})$   be a representation of  it.
Suppose $f:V\rightarrow D$ is a morphism between $D$-representations (from $V$ to the adjoint representation $D$), i.e.
 $$f(x\prec_{l}u)= x\prec  f(u), f(x\succ_{l}u)= x\succ  f(u), f(u\prec_{r}x)= f(u)\prec x,f( u\succ_{r}x)= f(u)\succ x$$
for all~$x,y\in D,u\in V. $
 Then there is a quadri-dendriform algebra structure
on $V$ with the operations $u\prec_{\vdash}v = f(u)\prec_l v,  u \succ_{\vdash}v= f(u)\succ_l v,  u\prec_{\dashv}v=u \prec_r  f(v) ~~\text{and}~~  u\succ_{\dashv} v=u\succ_r  f(v)$ for all $u,v\in V$.
\end{example}

Another example  emerges from the representations of a dendriform algebra.

\begin{proposition}\label{proposition:3-tri-Leibniz algebra}
Let $(V;\prec_{l},\succ_{l},\prec_{r},\succ_{r})$   be a  representation of a    dendriform algebra $(D,\prec,\succ)$.
Then $(D\oplus V, \prec_\vdash,\prec_\dashv,\succ_\vdash,\succ_\dashv)$ is a  quadri-dendriform algebra, where
 \begin{align*}
&(x,u)\prec_\vdash(y, v):=(x\prec y, x\prec_l v),~~~(x,u)\prec_\dashv(y, v):=(x\prec y, u\prec_r y),\\
&(x,u)\succ_\vdash(y, v):=(x\succ y, x\succ_l v),~~~(x,u)\succ_\dashv(y, v):=(x\succ y, u\succ_r y)
\end{align*}
for all $(x,u),(y,v)\in D\oplus V$.
This   quadri-dendriform algebra is called the hemisemidirect product  quadri-dendriform algebra and denoted by  $D\ltimes_{\mathrm{hemi}} V$.
\end{proposition}

\begin{proof}
For any $(x,u),(y,v),(z,w) \in  D\oplus V$ and $\ast\in \{ \vdash, \dashv\}$, by Eqs.  \eqref{2.1}-\eqref{2.12}, we have
\begin{align*}
& ((x,u)\prec_\ast(y,v))\prec_\vdash(z,w)=((x\prec y)\prec z,(x\prec y)\prec_l w)\\
&=(x\prec  (y\prec z+y\succ z),x\prec_{l}  (y\prec_{l} w+y\succ_{l} w))\\
&=(x,u)\prec_\vdash((y, v)\prec_\vdash(z,w)+(y, v)\succ_\vdash(z,w)).
\end{align*}
Similarly, we get
\begin{align*}
&((x,u)\succ_\ast  (y,v))\prec_{\vdash}  (z,w)=((x \succ y) \prec  z, (x \succ y)\prec_{l}w)\\
&=(x\succ  (y\prec z), x\succ_{l}  (y\prec_{l} w))=(x,u)\succ_{\vdash}  ((y,v)\prec_{\vdash} (z,w)),\\
&(x,u)\succ_{\vdash} ((y,v)\succ_{\vdash}  (z,w))=(x\succ(y\succ z),x\succ_{l} ( y \succ_{l}  w))  \\
&=((x\prec  y+x\succ  y)\succ z,(x\prec  y+x\succ  y)\succ_{l} w)  \\
&=\big((x,u)\prec_\ast  (y,v)+(x,u)\succ_\ast  (y,v)\big)\succ_{\vdash} (z,w),\\
&((x,u)\prec_{\vdash} (y,v))\prec_{\dashv}   (z,w)=( (x\prec y)\prec z, (x\prec_l v) \prec_{r}    z )\\
&=\big(x\prec  (y\prec z+y\succ z),x\prec_{l}  (v\prec_{r} z+v\succ_{r} z)\big)\\
&=(x,u)\prec_{\vdash}  \big((y,v)\prec_{\dashv}  (z,w)+(y,v)\succ_{\dashv}  (z,w)\big).
\end{align*}
Thus, the  Eqs.  \eqref{2.13}-\eqref{2.16} hold. The same for Eqs. \eqref{2.17}-\eqref{2.21}. This completes the proof.
\end{proof}

Now we introduce the notion of  relative  averaging
operators on  dendriform algebras, which are closely related to quadri-dendriform algebras.

\begin{definition}
Let $(D,\prec,\succ)$  be a  dendriform algebra  and  $(V;\prec_{l},\succ_{l},\prec_{r},\succ_{r})$   be a  representation  of  it.  A relative averaging
operator on $D$ with respect to the  representation  $(V;\prec_{l},\succ_{l},\prec_{r},\succ_{r})$
is a linear map
$T:V\rightarrow  D$ that satisfies
\begin{align}
&Tu\prec Tv =T (Tu\prec_{l}v)=T\rho_r(u\prec_{r}Tv), \label{2.22}\\
&Tu\succ Tv =T (Tu\succ_{l}v)=T\rho_r(u\succ_{r}Tv), \label{2.23}
\end{align}
for all  $u,v\in V$.
\end{definition}

The  relative averaging operator  on a  dendriform algebra   $(D,\prec,\succ)$   with respect to the adjoint  representation
$(D;\prec_{l}=\prec_{r}=\prec,\succ_{l}=\succ_{r}=\succ)$  is called an averaging operator. In this case the Eqs. \eqref{2.22} and \eqref{2.23}  can be written as
\begin{align*}
&Tx\prec Ty =T (Tx\prec y)=T(x\prec Ty),~~~Tx\succ Ty =T (Tx\succ y)=T(x\succ Ty),
\end{align*}
for all  $x,y\in D$.
 In \cite{Okba}, the authors  describe  all averaging operators on 2-dimensional dendriform algebras.

\begin{definition}
An averaging dendriform algebra  is a pair  $(D,T)$, consisting of a dendriform algebra     $(D,\prec,\succ)$   endowed with an
averaging operator $T$.
\end{definition}

\begin{example}
Let $(D,\prec,\succ)$  be a  dendriform algebra. Then the  map $k \mathrm{Id}: D\rightarrow D$ is an averaging
operator on $D$, where $\mathrm{Id}$ is an identity map and $k\in \mathbb{K}$.
\end{example}

\begin{example}
Let $(D,\prec,\succ,d)$  be a differential dendriform algebra  as presented in Example  \ref{exam:differential  dendriform algebra}.
Then    $d$ is an averaging operator on $D$.
\end{example}

\begin{example}
Let $(D,\prec,\succ)$  be a  dendriform algebra and  $(V;\prec_{l},\succ_{l},\prec_{r},\succ_{r})$   be a  representation  of  it.
Suppose $f:V\rightarrow D$ is a morphism of $D$-representations  from $V$ to the adjoint representation $D$.
 Then $f$ is a  relative averaging operator  on $D$ with respect to the representation $V$.
\end{example}

In the following, we give a characterization of a relative averaging operator  in terms
of the graph of the operator.

\begin{theorem} \label{theorem:graph}
A linear map $T:V\rightarrow  D$ is a relative averaging operator   on a dendriform algebra   $(D,\prec,\succ)$  with respect to the
representation   $(V;\prec_{l},\succ_{l},\prec_{r},\succ_{r})$    if and only if the graph $Gr(T)=\{(Tu,u)~|~u\in V\}$ is a subalgebra of the hemisemidirect product
quadri-dendriform algebra $D\ltimes_{\mathrm{hemi}} V$.
\end{theorem}

\begin{proof}
Let $T:V\rightarrow  D$  be a linear map.
For any $u,v\in V$,  we have
\begin{align*}
&(Tu,u)\prec_{\vdash}(Tv,v)=(Tu\prec Tv,Tu\prec_l v ),~~~(Tu,u)\prec_{\dashv}(Tv,v)=(Tu\prec Tv, u\prec_r Tv ),\\
&(Tu,u)\succ_{\vdash}(Tv,v)=(Tu\succ Tv,Tu\succ_l v ),~~~(Tu,u)\succ_{\dashv}(Tv,v)=(Tu\succ Tv, u\succ_r Tv ).
\end{align*}
The above four elements are in $Gr(R)$ if and only if $Tu\prec Tv=T(Tu\prec_l v)$, $Tu\prec Tv=T(u\prec_r Tv)$, $Tu\succ Tv=T(Tu\succ_l v)$  and $ Tu\succ Tv=T( u\succ_r Tv)$.
Therefore,  the graph $Gr(T)$ is a subalgebra of the hemisemidirect product  quadri-dendriform algebra $D\ltimes_{\mathrm{hemi}} V$ if and only if   $T$ is a relative averaging operator.
\end{proof}

Since  $Gr(T)$ is isomorphic to $V$ as vector spaces, so there is an induced
quadri-dendriform algebra structure on $V$.

\begin{proposition} \label{proposition:ET}
Let $T:V\rightarrow  D$ be a relative averaging operator  on a dendriform algebra   $(D,\prec,\succ)$  with respect to the
representation   $(V;\prec_{l},\succ_{l},\prec_{r},\succ_{r})$.  Then  $V$ inherits  a quadri-dendriform algebra structure with the operations
\begin{align*}
 &u \prec_{\vdash}^Tv=Tu\prec_l v,~~u\prec_{\dashv}^Tv=u\prec_r Tv,~~u\succ_{\vdash}^Tv=Tu\succ_l v,~~u\succ_{\dashv}^Tv= u\succ_r Tv,
\end{align*}
for  all $u,v\in V$.   Moreover, $T$ is a homomorphism from the quadri-dendriform algebra $(V, \prec_{\vdash}^T,\prec_{\dashv}^T,\succ_{\vdash}^T,\succ_{\dashv}^T)$ to the dendriform algebra  $(D,\prec,\succ)$.
\end{proposition}

\begin{cor} \label{cor:di-dendriform algebra}
Let $(D,\prec,\succ)$  be a dendriform algebra and $T: D\rightarrow  D$  be an averaging operator. Then,
$D$ is also equipped with  a quadri-dendriform algebra structure defined by
\begin{align*}
 &x \prec_{\vdash}^Ty=Tx\prec y,~~x\prec_{\dashv}^Ty=x\prec Ty,~~x\succ_{\vdash}^Ty=Tx\succ y,~~x\succ_{\dashv}^Ty=x\succ Ty,
\end{align*}for  all $x,y\in D$.
\end{cor}

The above result  indicates that  a relative averaging operator  on a dendriform algebra  with respect to a representation induces a  quadri-dendriform algebra  structure on the underlying representation  space.
 The subsequent result establishes the converse relationship.

\begin{theorem}  \label{theorem:di-dendriform algebrarelative averaging operator}
Every  quadri-dendriform algebra  can induce  a relative averaging operator  on a dendriform algebra
with respect to a representation.
\end{theorem}

\begin{proof}
Let  $(\mathcal{D},\prec_\vdash,\prec_\dashv,\succ_\vdash,\succ_\dashv)$   be a  quadri-dendriform algebra.
 Then by Example \ref{exam:3-tri-Leibniz kernel}, the quotient algebra $\mathcal{D}_{\mathrm{Dend}}=\frac{\mathcal{D}}{I_{\mathcal{D}}}$ carries a
dendriform algebra structure with the operations
$$ \overline{x}\prec\overline{y} =\overline{ x\prec_\vdash y}=\overline{ x\prec_\dashv y },~~\overline{x}\succ\overline{y} =\overline{ x\succ_\vdash y}=\overline{ x\succ_\dashv y }$$
for all $x,y\in \mathcal{D}$, here $\overline{x}$ denotes the class of an element $x\in \mathcal{D}$.
We define four bilinear maps
$\prec_l,\succ_l: \mathcal{D}_{\mathrm{Dend}}\times \mathcal{D}\rightarrow \mathcal{D}$ and
$\prec_r,\succ_r: \mathcal{D}\times \mathcal{D}_{\mathrm{Dend}}\rightarrow \mathcal{D}$
by
\begin{align*}
 & \overline{x}\prec_l y = x\prec_\vdash  y,~~  \overline{x}\succ_l y = x\succ_\vdash  y,~~y\prec_r \overline{x} =y\prec_\dashv x,~~ y\succ_r \overline{x} =y\succ_\dashv x,
\end{align*}
for all $\overline{x}\in \mathcal{D}_{\mathrm{Dend}}$ and $y\in \mathcal{D}$.
It is easy to verify that the maps $\prec_l,\succ_l,\prec_r,\succ_r$
define a   representation of the dendriform algebra  $\mathcal{D}_{\mathrm{Dend}}$ on the vector
space  $\mathcal{D}$ . Moreover, the quotient map $T:  \mathcal{D}\rightarrow  \mathcal{D}_{\mathrm{Dend}}, x\mapsto \overline{x}$ is a  relative averaging operator  as
\begin{align*}
 Tx\prec Ty = \overline{x}\prec \overline{y}  =
\left\{ \begin{array}{ll}
=\overline{x\prec_\vdash y}=\overline{\overline{x}\prec_l y}=T(Tx\prec_l y),\\
=\overline{x\prec_\dashv y}=\overline{x\prec_r\overline{y}}=T(x\prec_r Ty),
 \end{array}  \right.\\
  Tx\succ Ty = \overline{x}\succ \overline{y}  =
\left\{ \begin{array}{ll}
=\overline{x\succ_\vdash y}=\overline{\overline{x}\succ_l y}=T(Tx\succ_l y),\\
=\overline{x\succ_\dashv y}=\overline{x\succ_r\overline{y}}=T(x\succ_r Ty),
 \end{array}  \right.
 \end{align*}
 for all $x,y \in \mathcal{D}$.
\end{proof}

\begin{remark}
It is crucial to note that an arbitrary  quadri-dendriform algebra   may not be induced from an averaging   dendriform algebra.
However, any  quadri-dendriform algebra   can be embedded into an averaging  dendriform algebra.
More precisely, let   $(\mathcal{D},\prec_\vdash,\prec_\dashv,\succ_\vdash,\succ_\dashv)$   be a quadri-dendriform algebra.  Then by Proposition \ref{proposition:  dendriform algebra},
the direct sum $\mathcal{D}_{\mathrm{Dend}}\oplus \mathcal{D}$
inherits a  dendriform algebra structure with the operations
\begin{align*}
 (\overline{x},x)\prec_{\ltimes}(\overline{y},y) =( \overline{x}\prec\overline{y},\overline{x}\prec_l y+x\prec_r\overline{y}),\\
  (\overline{x},x)\succ_{\ltimes}(\overline{y},y) =( \overline{x}\succ\overline{y},\overline{x}\succ_l y+x\succ_r\overline{y}),
 \end{align*}
 for all $(\overline{x},x),(\overline{y},y)\in \mathcal{D}_{\mathrm{Dend}}\oplus \mathcal{D}.$
 Moreover, the map
 $$T:\mathcal{D}_{\mathrm{Dend}}\oplus \mathcal{D}\rightarrow \mathcal{D}_{\mathrm{Dend}}\oplus \mathcal{D},~~(\overline{x},y)\mapsto (\overline{y},0)$$
  is an averaging operator. Thus, $(\mathcal{D}_{\mathrm{Dend}}\oplus \mathcal{D},T)$ is an averaging dendriform algebra.
  Then the inclusion map
  $$i:\mathcal{D}\rightarrow \mathcal{D}_{\mathrm{Dend}}\oplus \mathcal{D},~~ x\mapsto(0,x)$$
  is an embedding of the quadri-dendriform algebra $(\mathcal{D},\prec_\vdash,\prec_\dashv,\succ_\vdash,\succ_\dashv)$ into the
averaging  dendriform algebra $(\mathcal{D}_{\mathrm{Dend}}\oplus \mathcal{D},T)$.
\end{remark}

\section{Six-dendriform algebras }\label{section:Six-dendriform algebras}
\def\theequation{\arabic{section}.\arabic{equation}}
\setcounter{equation} {0}
In this section,  we introduce the notion of a six-dendriform algebra, which is a splitting of a tri-associative algebra. Then,
 we show that  a homomorphic  relative averaging operator induces a six-dendriform algebra.

\begin{definition}
Let  $(D,\prec,\succ)$ and $(D',\prec',\succ')$  be two  dendriform algebras. We say that $D$ acts on $D'$
if there  four   bilinear maps:
\begin{align*}
\prec_l,\succ_l:D\otimes D'\rightarrow D',~~\prec_r,\succ_r:D'\otimes D\rightarrow D',
\end{align*}
such that $(D';\prec_l,\succ_l,\prec_r,\succ_r)$ is a  representation   of $D$ and for any
$x,y,z\in D$ and $u,v,w\in D'$,
\begin{align}
&(x\prec_lv)\prec'w=x\prec_l  \big( v\prec'w + v\succ'w\big),\label{3.1}\\
&(x\succ_lv )\prec'w=x\succ_{l}(v\prec'w),\label{3.2}\\
&x\succ_l ( v\succ'w)=(x\prec_lv+x\succ_lv)\succ' w,\label{3.3}\\
& (u\prec_ry)\prec'w=u\prec'\big(y\prec_lw+y\succ_lw\big),\label{3.4}\\
&(u\succ_ry)\prec'w=u\succ' (y\prec_lw),\label{3.5}\\
&u\succ' (y\succ_lw )=(u\prec_ry+ u\succ_ry)\succ' w,\label{3.6}\\
&(u\prec'v)\prec_r z=u\prec'\big( v\prec_rz+ v\succ_rz\big),\label{3.7}\\
&(u\succ'v)\prec_rz=u\succ' (v\prec_rz ),\label{3.8}\\
&u\succ' (v\succ_rz)=( u\prec'v+ u\succ'v)\succ_rz.\label{3.9}
\end{align}
\end{definition}

\begin{proposition}
Let  $(D,\prec,\succ)$ and $(D',\prec',\succ')$  be two  dendriform algebras,  and let
$\prec_l,\succ_l:D\otimes D'\rightarrow D'$ and $\prec_r,\succ_r:D'\otimes D\rightarrow D'$  be bilinear maps.
Then, the tuple $(D',\prec',\succ';\prec_l,\succ_l,\prec_r,\succ_r)$ is an action of $D$ if and only if
$D\oplus D'$ carries a dendriform algebra structure with operations given by
\begin{align*}
  (x,u)\prec_{\bowtie}(y,v) =( x\prec y, x\prec_lv+ u\prec_ry+ u\prec'v),\\
    (x,u)\succ_{\bowtie}(y,v) =( x\succ y, x\succ_lv+ u\succ_ry+ u\succ'v),
\end{align*}
for all $(x,u),(y,v) \in D\oplus D'$, which is called the semi-direct product of $D$ with $D'$.
\end{proposition}

\begin{proof}
Let  $(D,\prec,\succ)$ and $(D',\prec',\succ')$  be two  dendriform algebras.
Then, for any $(x,u),(y,v), (z,w) \in D\oplus D'$, we have
\begin{align*}
& \big((x,u)\prec_{\bowtie} (y,v)\big)\prec_{\bowtie}  (z,w)-(x,u)\prec_{\bowtie}  \big((y,v)\prec_{\bowtie} (z,w) +(y,v)\succ_{\bowtie} (z,w) \big)\\
&= \Big( (x\prec y)\prec z-x\prec(y\prec z+ y\succ z),~~ (x\prec y)\prec_l w+(x\prec_lv+ u\prec_ry+ u\prec'v)\prec_r z+\\
 &~~~~(x\prec_lv+ u\prec_ry+ u\prec'v)\prec'w - x\prec_l  \big(y\prec_lw+ v\prec_rz+ v\prec'w +y\succ_lw+ v\succ_rz+ v\succ'w\big)-\\
&~~~~~u\prec_r(y\prec z+ y\succ z)-u\prec'\big(y\prec_lw+ v\prec_rz+ v\prec'w +y\succ_lw+ v\succ_rz+ v\succ'w\big)\Big)\\
&= \Big(0,~(x\prec y)\prec_l w+(x\prec_lv+ u\prec_ry+ u\prec'v)\prec_r z+(x\prec_lv+ u\prec_ry)\prec'w - \\
&~~~~~x\prec_l  \big(y\prec_lw+ v\prec_rz+ v\prec'w +y\succ_lw+ v\succ_rz+ v\succ'w\big)-u\prec_r(y\prec z+ y\succ z)-\\
&~~~~~u\prec'\big(y\prec_lw+ v\prec_rz +y\succ_lw+ v\succ_rz\big)\Big).\\
\end{align*}
Similarly, we obtain
\begin{align*}
&\big((x,u)\succ_{\bowtie}  (y,v)\big)\prec_{\bowtie}  (z,w) -(x,u)\succ_{\bowtie}  \big((y,v)\prec_{\bowtie} (z,w) \big)\\
&=\Big((x\succ y)\prec z-x\succ(y\prec z),~~(x\succ y)\prec_lw+\big(x\succ_lv+ u\succ_ry+ u\succ'v\big)\prec_rz+\\
&~~~~~\big(x\succ_lv+ u\succ_ry+ u\succ'v\big)\prec'w-x\succ_{l} \big(y\prec_lw+ v\prec_rz+ v\prec'w \big)-u\succ_r(y\prec z)-\\
&~~~~~u\succ'\big(y\prec_lw+ v\prec_rz+ v\prec'w \big)\Big)\\
&=\Big(0,~(x\succ y)\prec_lw+\big(x\succ_lv+ u\succ_ry+ u\succ'v\big)\prec_rz+\big(x\succ_lv+ u\succ_ry\big)\prec'w-\\
&~~~~~x\succ_{l} \big(y\prec_lw+ v\prec_rz+ v\prec'w \big)-u\succ_r(y\prec z)-u\succ'\big(y\prec_lw+ v\prec_rz\big)\Big)
\end{align*}
and
\begin{align*}
&(x,u)\succ_{\bowtie} \big((y,v)\succ_{\bowtie}  (z,w)\big)-\big((x,u)\prec_{\bowtie}  (y,v)+(x,u)\succ_{\bowtie}  (y,v)\big)\succ_{\bowtie} (z,w)\\
&=\Big(x\succ( y\succ z)-\big(x\prec y+x\succ y\big)\succ z,~~x\succ_l\big(y\succ_lw+ v\succ_rz+ v\succ'w\big)+u\succ_r( y\succ z)+\\
&~~~~u\succ'\big(y\succ_lw+ v\succ_rz+ v\succ'w\big)-(x\prec y)\succ_lw-\big(x\prec_lv+ u\prec_ry+ u\prec'v\big)\succ_rz-\\
&~~~~\big(x\prec_lv+ u\prec_ry+ u\prec'v\big)\succ' w-(x\succ y)\succ_l w-\big(x\succ_lv+ u\succ_ry+ u\succ'v\big)\succ_{r} z-\\
&~~~~\big(x\succ_lv+ u\succ_ry+ u\succ'v\big)\succ' w\Big)\\
&=\Big(0,~x\succ_l\big(y\succ_lw+ v\succ_rz+ v\succ'w\big)+u\succ_r( y\succ z)+u\succ'\big(y\succ_lw+ v\succ_rz\big)-\\
&~~~~(x\prec y)\succ_lw-\big(x\prec_lv+ u\prec_ry+ u\prec'v\big)\succ_rz-\big(x\prec_lv+ u\prec_ry\big)\succ' w-\\
&~~~~(x\succ y)\succ_l w-\big(x\succ_lv+ u\succ_ry+ u\succ'v\big)\succ_{r} z-\big(x\succ_lv+ u\succ_ry\big)\succ' w\Big).
\end{align*}
Then, $(D\oplus D',  \prec_{\bowtie},\succ_{\bowtie})$ is a dendriform algebra if and only if  $(D';\prec_l,\succ_l,\prec_r,\succ_r)$ is a  representation   of $D$ and
Eqs. \eqref{3.1}-\eqref{3.9} are satisfied. This completes the proof.
\end{proof}

\begin{definition}
 A linear map $T:D'\rightarrow D$ is called a homomorphic  relative averaging operator  on the dendriform algebra $(D,\prec,\succ)$    with respect to a   action  $(D',\prec',\succ';\prec_l,\succ_l,\prec_r,\succ_r)$    if $T$ qualifies both as a relative averaging operator and a  dendriform  homomorphism.
\end{definition}

A homomorphic  relative averaging operator on a   dendriform algebra is also called  a relative averaging operator of  weight $1$ on a   dendriform algebra.  In particular, a homomorphic relative averaging operator  on a
dendriform algebra $(D,\prec,\succ)$   with respect to the adjoint  representation
is called a homomorphic averaging operator.

\begin{example}
Let  $(D,\prec,\succ)$    be a  dendriform algebra. Then the  identity map $ \mathrm{Id}: D\rightarrow D$ is a  homomorphic averaging
operator on $D$.
\end{example}

\begin{example}
Let $(D,\prec,\succ)$    be a  dendriform algebra  and   $(D',\prec',\succ';\prec_l,\succ_l,\prec_r,\succ_r)$   be an action of  $D$.
 Suppose $f:D'\rightarrow D$
is a dendriform  algebra homomorphism and a morphism of $D$-representations.  Then $f$ is a  homomorphic relative averaging operator.
\end{example}

\begin{example}
Let $(D,\prec,\succ)$    be a  dendriform algebra.  Then the dendriform algebra operations of $D$ induces a
dendriform operation on the space $F=D[\![t]\!]/(t^2)$ using $t$-bilinearity. Further, $F$ can be given an
action of  $D$ with the  four bilinear maps are respectively given by
\begin{align*}
&x\prec_l(a+tb)=x\prec a+t(x\prec b),~~x\succ_l(a+tb)=x\succ a+t(x\succ b),\\
&(a+tb)\prec_r x=a\prec x+t(b\prec x),~~(a+tb)\succ_r x=a\succ   x+t(b\succ x),
\end{align*}
for all $x\in D$ and  $a+tb\in F$. With this, the projection map $T:F\rightarrow D$ given by $T(a+tb)=a$ is a homomorphic relative averaging operator.
\end{example}

 In \cite{Loday2004},   Loday and Ronco introduced the notion of   tri-associative algebras.
Inspired by Loday and Ronco's work, we   present the notion of six-dendriform algebras, which are closely related to  tri-associative algebras.

\begin{definition}
A   six-dendriform algebra is a tuple $(\mathcal{D},\prec_\perp,\succ_\perp,\prec_\vdash,\succ_\vdash,\prec_\dashv,\succ_\dashv)$ consisting of a dendriform algebra  $(\mathcal{D},\prec_\perp,\succ_\perp)$
 and a  quadri-dendriform algebra  $(\mathcal{D}, \prec_\vdash,\succ_\vdash,\prec_\dashv,\succ_\dashv)$   such that
 for all $x,y,z \in D$,
\begin{align}
&\left\{ \begin{array}{lllllllll}
(x\prec_\vdash y)\prec_\perp z=x\prec_\vdash  \big( y\prec_\perp z + y\succ_\perp z\big),\\
(x\succ_\vdash y )\prec_\perp z=x\succ_{\vdash}(y\prec_\perp z), \\
x\succ_\vdash ( y\succ_\perp z)=(x\prec_\vdash y+x\succ_\vdash y)\succ_\perp z, \\
 (x\prec_\dashv y)\prec_\perp z=x\prec_\perp\big(y\prec_\vdash z+y\succ_\vdash z\big), \\
(x\succ_\dashv y)\prec_\perp z=x\succ_\perp (y\prec_\vdash z), \\
x\succ_\perp (y\succ_\vdash z )=(x\prec_\dashv y+ x\succ_\dashv y)\succ_\perp z, \\
(x\prec_\perp y)\prec_\dashv z=x\prec_\perp\big( y\prec_\dashv z+ y\succ_\dashv z\big), \\
(x\succ_\perp y)\prec_\dashv z=x\succ_\perp (y\prec_\dashv z ), \\
x\succ_\perp (y\succ_\dashv z)=(x\prec_\perp y+ x\succ_\perp y)\succ_\dashv z,
 \end{array}  \right.\label{3.10}\\
 &\left\{ \begin{array}{llllllll}
(x\prec_\perp y)\prec_\vdash z=(x\prec_\vdash y)\prec_\vdash z=(x\prec_\dashv y)\prec_\vdash z,\\
(x\succ_\perp y)\prec_\vdash z=(x\succ_\vdash y)\prec_\vdash z=(x\succ_\dashv y)\prec_\vdash z,\\
(x\prec_\perp y)\succ_\vdash z=(x\prec_\vdash y)\succ_\vdash z=(x\prec_\dashv y)\succ_\vdash z,\\
(x\succ_\perp y)\succ_\vdash z=(x\succ_\vdash y)\succ_\vdash z=(x\succ_\dashv y)\succ_\vdash z,
 \end{array}  \right.\label{3.11}\\
  &\left\{ \begin{array}{llllllll}
x\prec_\dashv (y\prec_\perp z)=x\prec_\dashv (y\prec_\vdash z)=x\prec_\dashv (y\prec_\dashv z),\\
x\succ_\dashv (y\prec_\perp z)=x\succ_\dashv (y\prec_\vdash z)=x\succ_\dashv (y\prec_\dashv z),\\
x\prec_\dashv (y\succ_\perp z)=x\prec_\dashv (y\succ_\vdash z)=x\prec_\dashv (y\succ_\dashv z),\\
x\succ_\dashv (y\succ_\perp z)=x\succ_\dashv (y\succ_\vdash z)=x\succ_\dashv (y\succ_\dashv z).
 \end{array}  \right.\label{3.12}
\end{align}
\end{definition}
A morphism between two six-dendriform algebras is a linear map which
is compatible with the six operations. We denote by $\mathbf{SDend}$ the category
of six-dendriform algebras.

If $\prec_\perp=\prec_\vdash=\prec_\dashv$ and $\succ_\perp=\succ_\vdash=\succ_\dashv$, then we get simply a dendriform algebra.
So we get a functor between the categories of algebras :
$$\mathbf{Dend}\rightarrow\mathbf{SDend}.$$
Ignoring the operations $\prec_\perp,\succ_\perp$ gives a quadri-dendriform algebra. Hence there is a
(forgetful) functor
$$\mathbf{SDend}\rightarrow\mathbf{QDend}$$
from the category of six-dendriform algebras to the category of quadri-dendriform algebras.


\begin{definition}\cite{Loday2004}
A tri-associative algebra is a vector space $A$ equipped  with three   binary operations  $\vdash,\dashv,\perp: A\otimes A\rightarrow A$,  such that $(A,  \vdash, \dashv)$  constitutes a  di-associative algebra,
$(A,  \perp)$  forms an   associative algebra, and
for all $x, y, z\in A$,   the following equations hold:
\begin{align*}
& 
(x\dashv y)\dashv  z=x\dashv  (y\perp z), ~~~(x\perp  y)\dashv  z=x\perp  (y\dashv z),  ~~~(x\dashv  y)\perp  z=x\perp (y\vdash z),  \\
&(x\vdash  y)\perp z=x\vdash  (y\perp z),  ~~~(x\perp y)\vdash  z=x\vdash  (y\vdash z).
\end{align*}
\end{definition}

In \cite{Das2023,Qiao2023}, the authors   showed that tri-associative algebras arise  from (relative) averaging operators of any nonzero
weight on associative algebras.

\begin{proposition}\label{proposition: tri-associative algebra}
Let $(\mathcal{D},\prec_\perp,\succ_\perp,\prec_\vdash,\prec_\dashv,\succ_\vdash,\succ_\dashv)$   be a  six-dendriform algebra.
Then $(\mathcal{D}, \perp,\vdash, \dashv)$  is a  tri-associative algebra, where
 \begin{eqnarray*}
&&x \perp y=x\prec_\perp y+x\succ_\perp y,~~x \vdash y=x\prec_\vdash y+x\succ_\vdash y,~~x\dashv y=x \prec_\dashv y+x\succ_\dashv y,
\end{eqnarray*}
for all $x,y\in  \mathcal{D}$.
\end{proposition}

\begin{proof}
It follows from Proposition \ref{proposition: di-associative algebra} that $(\mathcal{D},  \vdash, \dashv)$  is a  di-associative algebra. 
The fact that $(\mathcal{D},\prec_\perp,\succ_\perp)$ is a dendriform algebra means that $(\mathcal{D}, \perp)$ is a   associative algebra.
Next, for any $x,y,z\in  \mathcal{D}$, by Eqs.  \eqref{2.13}-\eqref{2.21}, \eqref{3.10}-\eqref{3.12}, we have
\begin{align*}
& (x\dashv y)\dashv  z=(x \prec_\dashv y+x\succ_\dashv y)\prec_\dashv  z+(x \prec_\dashv y+x\succ_\dashv y)\succ_\dashv  z\\
&=x\prec_\dashv  (y\prec_\vdash z+y\succ_\vdash z)+x\succ_\dashv  (y\prec_\dashv z)+x\succ_\dashv  (y\succ_\vdash z)\\
&=x\prec_\dashv  (y\prec_\perp z+y\succ_\perp z)+x\succ_\dashv  (y\prec_\perp z+y\succ_\perp z)\\
&=x\dashv  (y\perp z).
\end{align*}
Similarly, we get
\begin{align*}
&(x\perp  y)\dashv  z=(x\prec_\perp y+x\succ_\perp y)\prec_\dashv  z+(x\prec_\perp y+x\succ_\perp y)\succ_\dashv  z\\
&=x\prec_\perp\big( y\prec_\dashv z+ y\succ_\dashv z\big)+x\succ_\perp (y\prec_\dashv z )+x\succ_\perp (y\succ_\dashv z)\\
&=x\perp  (y\dashv z),\\
&(x\dashv  y)\perp  z=(x \prec_\dashv y+x\succ_\dashv y)\prec_\perp  z+(x \prec_\dashv y+x\succ_\dashv y)\succ_\perp  z\\
&=x\prec_\perp (y\prec_\vdash z+y\succ_\vdash z)+x\succ_\perp (y\prec_\vdash z)+x\succ_\perp (y\succ_\vdash z)\\
&=x\perp (y\vdash z),\\
&(x\vdash  y)\vdash  z=(x\prec_\vdash y+x\succ_\vdash y)\prec_\vdash  z+(x\prec_\vdash y+x\succ_\vdash y)\succ_\vdash  z\\
&=x\prec_{\vdash}  (y\prec_{\vdash} z+y\succ_{\vdash} z)+x\succ_{\vdash}  (y\prec_{\vdash} z)+x\succ_{\vdash} (y\succ_{\vdash}  z)\\
&=x\vdash  (y\vdash z).
\end{align*}
Finally,  we   obtain
\begin{align*}
&(x\perp y)\vdash  z=(x\prec_\perp y+x\succ_\perp y)\prec_\vdash  z+(x\prec_\perp y+x\succ_\perp y)\succ_\vdash  z\\
&=(x\prec_\vdash y+x\succ_\vdash y)\prec_\vdash  z+(x\prec_\vdash y+x\succ_\vdash y)\succ_\vdash  z\\
&=x\prec_\vdash  (y\prec_\vdash z+y\succ_\vdash z)+x\succ_\vdash  (y\prec_\vdash z)+x\succ_\vdash  (y\succ_\vdash z)\\
&=x\vdash  (y\vdash z).
\end{align*}
This completes the proof.
\end{proof}
In other words, a six-dendriform algebra is a tri-associative algebra for
which the tri-associative operations are the sum of two operations.
We denote by $\mathbf{Triass}$ the category of tri-associative algebras. By the preceding proposition, there is a
well-defined functor:
$$\mathbf{SDend}\rightarrow\mathbf{Triass}.$$

\begin{theorem}
Let $T:D'\rightarrow D$ be a homomorphic  relative averaging operator on the
 dendriform algebra $(D,\prec,\succ)$    with respect to a   action  $(D',\prec',\succ';\prec_l,\succ_l,\prec_r,\succ_r)$.
Then $(D',\prec_\perp,\succ_\perp, \prec_\vdash^T,\prec_\dashv^T,\succ_\vdash^T,\succ_\dashv^T)$ is a six-dendriform algebra, where
\begin{align*}
&u \prec_\perp v= u\prec' v,~~u \succ_\perp v= u\succ' v,~~u \prec_{\vdash}^Tv=Tu\prec_l v,\\
&u\prec_{\dashv}^Tv=u\prec_r Tv,~~u\succ_{\vdash}^Tv=Tu\succ_l v,~~u\succ_{\dashv}^Tv= u\succ_r Tv,
\end{align*}
for all $u,v \in D'$.
\end{theorem}

\begin{proof}
We have already $(D',\prec_\perp,\succ_\perp)$ is a dendriform algebra. Moreover, it follows from Proposition \ref{proposition:ET} that
$(D', \prec_\vdash^T,\prec_\dashv^T,\succ_\vdash^T,\succ_\dashv^T)$ is a quadri-dendriform algebra. Now, for any  $u,v,w \in D'$, according to Eqs.  \eqref{3.1}-\eqref{3.9}, we have
\begin{align*}
&(u\prec_\vdash^T v)\prec_\perp w=(Tu\prec_l v)\prec' w=Tu\prec_l  \big( v\prec' w + v\succ' w\big)=u\prec_\vdash^T  \big( v\prec_\perp w + v\succ_\perp w\big),\\
&(u\succ_\vdash^T v )\prec_\perp w=(Tu\succ_l  v )\prec' w=Tu\succ_{l}(v\prec'w)=u\succ_{\vdash}^T(v\prec\perp w), \\
&u\succ_\vdash^T ( v\succ_\perp w)=Tu\succ_l  ( v\succ' w)=(Tu\prec_lv+Tu\succ_lv)\succ' w=(u\prec_\vdash^T v+u\succ_\vdash^T v)\succ_\perp w, \\
&(u\prec_\dashv^T v)\prec_\perp w=(u\prec_r Tv)\prec' w=u\prec'\big(Tv\prec_lw+Tv\succ_lw\big)=u\prec_\perp\big(v\prec_\vdash^T w+v\succ_\vdash^T w\big), \\
&(u\succ_\dashv^T v)\prec_\perp w=(u\succ_r Tv)\prec' w=u\succ' (Tv\prec_lw)=u\succ_\perp (v\prec_\vdash^T w), \\
&u\succ_\perp (v\succ_\vdash^T w )=u\succ' (Tv\succ_l w )=(u\prec_rTv+ u\succ_rTv)\succ' w=(u\prec_\dashv^T v+ u\succ_\dashv^T v)\succ_\perp w, \\
&(u\prec_\perp v)\prec_\dashv^T w=(u\prec' v)\prec_r T w=u\prec'\big( v\prec_rT w+ v\succ_rT w\big)=u\prec_\perp\big( v\prec_\dashv^T w+ v\succ_\dashv^T w\big), \\
&(u\succ_\perp v)\prec_\dashv^T w=(u\succ' v)\prec_r T w= u\succ' (v\prec_r Tw )=u\succ_\perp (v\prec_\dashv^T w ), \\
&u\succ_\perp (v\succ_\dashv^T w)=u\succ'(v\succ_r T w)=( u\prec'v+ u\succ'v)\succ_r Tw=(u\prec_\perp v+ u\succ_\perp v)\succ_\dashv^T w.
\end{align*}
On the other hand, we have
 \begin{align*}
& (u\prec_\perp v) \prec_\vdash^T w=T(u\prec' v)\prec_\vdash w =(Tu\prec T v)\prec_\vdash w=
\left\{ \begin{array}{ll}
=T(Tu\prec_l v)\prec_l  w=(u\prec_\vdash^T v) \prec_\vdash^T w,\\
=T(u\prec_r Tv)\prec_l  w=(u\prec_\dashv^T v) \prec_\vdash^T w,
 \end{array}  \right.
\end{align*}
Therefore,  Eqs. \eqref{3.10}  and the first equation in Eqs.  \eqref{3.11}  hold.  Similarly, it can be proved that other equations in Eqs.  \eqref{3.11}  and Eqs.  \eqref{3.12} are also true.
  This completes the proof.
\end{proof}

\begin{cor} \label{cor:six-dendriform algebra}
Let $(D,\prec,\succ)$  be a dendriform algebra and $T: D\rightarrow  D$  be a homomorphic averaging operator. Then,
$(D,\prec_\perp,\succ_\perp, \prec_\vdash^T,\prec_\dashv^T,\succ_\vdash^T,\succ_\dashv^T)$ is a six-dendriform algebra, where
\begin{align*}
&\prec_\perp = \prec,~\succ_\perp = \succ,~x \prec_{\vdash}^Ty=Tx\prec y,~x\prec_{\dashv}^Ty=x\prec Ty,~x\succ_{\vdash}^Ty=Tx\succ y,~x\succ_{\dashv}^Ty= x\succ Ty,
\end{align*}
for all  $x,y\in D$.
\end{cor}

The above result demonstrates that a homomorphic relative averaging operator on a dendriform algebra
with respect to a action induces a six-dendriform algebra structure on the underlying action space.

Conversely,  every six-dendriform algebra can induce a homomorphic relative averaging operator
on a dendriform algebra with respect to a action.
Let $(\mathcal{D},\prec_\perp,\succ_\perp,\prec_\vdash,\succ_\vdash,\prec_\dashv,\succ_\dashv)$ be a six-dendriform algebra. Let $\mathcal{I}$ be the ideal of $\mathcal{D}$ generated by the set $\big\{x\prec_\vdash y-x\prec_\dashv y,~~x\succ_\vdash y-x\succ_\dashv y~:~x,y\in \mathcal{D}\big\}$. Then the quotient space $\frac{\mathcal{D}}{\mathcal{I}}$ has the   dendriform algebra   algebra structure with the
operations
$$ \overline{x}\prec\overline{y} =\overline{ x\prec_\vdash y}=\overline{ x\prec_\dashv y },~~\overline{x}\succ\overline{y} =\overline{ x\succ_\vdash y}=\overline{ x\succ_\dashv y },$$
for all $ \overline{x},\overline{y}\in  \frac{\mathcal{D}}{\mathcal{I}}$.
We denote this dendriform algebra simply by $\mathcal{D}_{\mathrm{Dend}}$. On the other hand, the vector space $\mathcal{D}$ has the
dendriform algebra operations $\prec_\perp,\succ_\perp$. Denote this dendriform algebra $\mathcal{D}_\perp$. Moreover, the dendriform algebra $\mathcal{D}_\perp$ can be
given an action of $\mathcal{D}_{\mathrm{Dend}}$ given by
$$\overline{x}\prec_l y=x\prec_\vdash y,\overline{x}\succ_ly=x\succ_\vdash y,~y\prec_r\overline{x}=y\prec_\dashv x,y\succ_r\overline{x}=y\succ_\dashv x,$$
With these notations, the quotient map $T:\mathcal{D}_\perp\rightarrow \mathcal{D}_{\mathrm{Dend}}, x\mapsto \overline{x}$, for $x\in \mathcal{D}_\perp$, is a homomorphic  relative averaging operator.



\begin{center}
 {\bf ACKNOWLEDGEMENT}
 \end{center}

  The paper is  supported by the National Natural Science Foundation of China (Grant No. 12361005) and the  Universities Key Laboratory of System Modeling and Data Mining in Guizhou Province (Grant No.\, 2023013).

\renewcommand{\refname}{REFERENCES}

\end{document}